\documentclass[12pt]{amsart}
\usepackage{amscd}

\def\NZQ{\Bbb}               
\def\NN{{\NZQ N}}

\def\PP{{\Bbb P}}
\def\ZZ{{\NZQ Z}}
\def\RR{{\NZQ R}}
\def\D{{\Delta}}

\def\a{{\bold a}}
\def\b{{\bold b}}
\def\c{{\bold c}}
\def\w{{\bold w}}
\def\x{{\bold x}}
\def\1{{\mathbf 1}}
\def\0{{\mathbf 0}}

\newtheorem{Theorem}{Theorem}[section]
\newtheorem{Lemma}[Theorem]{Lemma}
\newtheorem{Corollary}[Theorem]{Corollary}

\newtheorem{Example}[Theorem]{Example}

\begin{document}

\title{Vertex cover algebras \\ of unimodular hypergraphs}

\author{J\"urgen Herzog, Takayuki Hibi and  Ng\^o Vi\^et Trung}

\address{J\"urgen Herzog, Fachbereich Mathematik und
Informatik, Universit\"at Duisburg-Essen, Campus Essen, 45117
Essen, Germany} \email{juergen.herzog@uni-essen.de}

\address{Takayuki Hibi, Department of Pure and Applied Mathematics,
Graduate School of Information Science and Technology, Osaka
University, Toyonaka, Osaka 560-0043, Japan}
\email{hibi@math.sci.osaka-u.ac.jp}

\address{Ng\^o Vi\^et Trung, Institute of Mathematics,
Vien Toan Hoc, 18 Hoang Quoc Viet, 10307 Hanoi,
Vietnam} \email{nvtrung@math.ac.vn}

 \maketitle

\begin{abstract}
It is proved that all vertex cover algebras of a hypergraph are standard graded if and only if the hypergraph is unimodular. This has interesting consequences on the symbolic powers of monomial ideals.
\end{abstract}
\bigskip

\section{Introduction}

A hypergraph $\D$ is a collection of subsets of a finite set of vertices $V$.
These subsets are called the edges of $\D$. For convenience, we assume throughout this paper that there is no inclusion between the edges of $\D$. Such a hypergraph is also called a clutter.

A vertex cover of $\D$ is a subset of $V$ which meets every edge of $\D$.
Suppose that $V = \{1,...,n\}$. We may think of a vertex cover of $V$ as a (0,1) vector $\c = (c_1,...,c_n)$ that satisfies the condition $\sum_{i \in F}c_i \ge 1$ for all $F \in \D$. 

Let $w\!: F \mapsto w_F $ be a weight function from $\D$ to the set of positive integers. We call $(\D,w)$ a weighted hypergraph. For $k \in \NN$ we define a $k$-cover of $(\D,w)$ as a vector $\c   \in \NN^n$ that satisfies the condition $\sum_{i \in F}c_i \ge kw_F$ for all $F \in \D$. 

Let $R = K[x_1,...,x_n]$ be a polynomial ring over a field $K$. For very vector $\c   \in \NN^n$ we set $\x^\c   = x_1^{c_1}\cdots x_n^{c_n}$. The {\it vertex cover algebra} $A(\D,w)$ is defined as the subalgebra of the one variable polynomial ring $R[t]$ generated by all monomials $\x^\c  t^k$, where $\c  $ is a $k$-cover of $(\D,w)$, $k \ge 0$. This algebra was introduced in \cite{HHT}. It was proved there that $A(\D,w)$ is a finitely generated normal Cohen-Macaulay ring \cite[Theorem 4.2]{HHT}.

For every edge $F \in \D$ let $P_F$ denote the ideal of $R$ generated by the variables $x_i$, $i \in F$. Let 
$I(\D,w) := \bigcap_{F \in \D}P_F^{w_F}.$
It is not hard to see that $A(\D,w)$ is the symbolic Rees algebra of the ideal
$I(\D,w)$, that is,
$$A(\D,w) = \oplus_{k \ge 0}I(\D,w)^{(k)}t^k,$$
where $I(\D,w)^{(k)}$ denotes the $k$-th symbolic power $\bigcap_{F \in \D}P_F^{w_Fk}$ of $I(\D,w)$.

We may view $A(\D,w)$ as a graded algebra over $R$ with $A(\D,w)_k =  I(\D,w)^{(k)}t^k$ for $k \ge 0$. It is obvious that  $A(\D,w)$ is  standard graded ($A(\D,w)$ is generated by elements of degree one) if and only if 
$I(\D,w)^{(k)} = I(\D,w)^k$ for all $k \ge 0$. Therefore, it is of great interest to know when $A(\D,w)$ is a standard graded algebra over $R$.

This problem has a satisfactory answer in the case $w$ is the canonical weight, i.e. $w_F = 1$ for all $F \in \D$.
In this case, we will use the notation $A(\D)$ and $I^*(\D)$ instead of $A(\D,w)$ and $I(\D,w)$. Notice that
ideals of the form $I^*(\D)$ are exactly squarefree monomial ideals of $R$. 
It was proved in \cite[Theorem 1.4]{HHTZ} (and implicitly in  \cite[Proposition 3.4]{Es} and \cite[Theorem 3.5]{Gi}) that $A(\D)$ is standard graded over $R$ if and only if the blocker of $\D$ is a Mengerian hypergraph (see the last section for more details). 

The aim of this paper is to characterize hypergraphs $\D$ for which $A(\D,w)$ is a standard graded algebra over $R$ for all weight functions $w$. If $\D$ is a graph, such a characterization was already obtained in \cite[Theorem 5.1 and Theorem 5.4]{HHT}, where it is proved that $A(\D,w)$ is a standard graded algebra over $R$ for all weight functions $w$ if and only if $\D$ is bipartite.
This result can be fully  generalized as follows.

Let $\D = \{F_1,...,F_m\}$. Let 
$M$ denote the edge-vertex incidence matrix of $\D$, that is,
 $M = (e_{ij})$, $i = 1,...,m$, $j = 1,...,n$, with $e_{ij} =  1$ if $i \in F_j$ and $e_{ij} = 0$ if $i \not\in F_j$. 
We say that $\D$ is an {\it unimodular hypergraph} if $M$ is totally unimodular, i.e. 
each subdeterminant of $M$ is 0,$\pm 1$. 

\begin{Theorem} \label{main}
 $A(\D,w)$ is a standard graded algebra for all weight functions $w$ if and only if 
 $\D$ is an unimodular hypergraph.
\end{Theorem}

It is known that if $\D$ has no alternating chain $v_1,F_1,v_2,F_2,...,v_s,F_s,x_{s+1}=x_1$ of odd length $s \ge 3$, where $x_1,...,x_s$ and $F_1,...,F_s$ are different vertices and edges of $\D$ and $x_i,x_{i+1} \in F_i$ for all $i = 1,...,s$, then $\D$ is unimodular (see e.g. \cite[Theorem 5, p.164]{Be}). In particular,  a graph is unimodular if and only if it is bipartite. Hence the afore mentioned result of \cite{HHT} is a special case of Theorem \ref{main}. 

An interesting case of Theorem \ref{main} concerns the hypergraph of all $n$-subsets of $n+1$ vertices.
In this case, $I^*(\D)$ is the intersection of the defining ideals of the $n+1$ points $(0,...,1,...,0)$ in $\PP^n$. 
Since every set of $n+1$ points in general position in $\PP^n$ can be transformed into this case, we obtain the following consequence.

\begin{Corollary}
Let $P_0,...,P_n$ be the defining ideals of $n+1$ points in general position in ${\PP}^n$. Then
$$(P_0^{w_0} \cap \cdots \cap P_n^{w_n})^k = P_0^{w_0k} \cap \cdots \cap P_n^{w_nk}$$
for all integers $w_0,...,w_s \ge 1$ and $k \ge 0$.
\end{Corollary}

\section{Decomposition property versus unimodularity}

We adhere to the notions of the preceding section.

Let $\D = \{F_1,...,F_m\}$ be a hypergraph on the set of vertices $\{1,...,n\}$. Let $w$ be a weight function on $\D$ and $\w = (w_{F_1},...,w_{F_m})$. Note that $\w \in \ZZ_+^m$. Let $M$ be the edge-vertex incidence matrix of $\D$. By definition, a vector $\c   \in \NN^n$ is a $k$-cover of the weighted hypergraph $(\D,w)$ if and only if $M\cdot \c   \ge k\w$. 

It is obvious  that the vertex cover algebra $A(\D,w)$ is standard graded over $R$ if and only if every monomial of $A(\D,w)_k$ is the product of $k$ monomials of $A(\D,w)_1$ for all $k \ge 1$ or, equivalently, if and only if
every $k$-cover of $(\D,w)$ can be written as a sum of $k$ 1-covers of $(\D,w)$ for all $k \ge 1$. This observation leads us to the following notion in hypergraph theory.

Let $Q$ be an arbitrary polyhedron in $\RR^n$ and $kQ = \{k\c|\ \c \in Q\}$, $k \ge 1$.
We say that $Q$ has the {\it integer decomposition property} if for each $k \ge 1$ and each integral vector $\c   \in kQ$ there exist integral vectors $\c_1,...,\c_k \in Q$ such that $\c   = \c_1 + \cdots +\c_k$. 

If $Q = \{\c   \in \RR^n|\ M\cdot \c   \ge \w \}$, then $kQ = \{\c   \in \RR^n|\ M\cdot \c   \ge k\w \}.$ Therefore, we have the following lemma.

\begin{Lemma} \label{decomposition}
$A(\D,w)$ is standard graded over $R$ if and only if the polyhedron  $\{\c   \in \RR^n|\ M\cdot \c   \ge \w \}$ has the integer decomposition property.
\end{Lemma}

If $Q = \{\c   \in \RR^n|\ M\cdot \c   \le \w \}$, there is the following characterization of the integer decomposition property in terms of $M$.

\begin{Theorem} \label{BT} {\rm (Baum and Trotter, see e.g. \cite[Theorem 19.4]{Sc})}
Let $M$ be an $m \times n$ integral matrix. 
Then the polyhedron $\{\c   \in \RR^n|\ M\cdot \c   \le \w \}$ has the decomposition property for all integral vector $\w \in \ZZ^m$ if and only if $M$  is totally unimodular.
\end{Theorem}

Since $M$ is totally unimodular if and only if $-M$ is totally unimodular and since $\w  $ may have negative components, we can replace $M$ by $-M$ and $\w  $ by $-\w  $ in Theorem \ref{BT}.  Therefore, Theorem \ref{BT} remains true if we replace the the polyhedron $\{\c   \in \RR^n|\ M\cdot \c   \le \w \}$ by the polyhedron $\{\c   \in \RR^n|\ M\cdot \c   \ge \w \}$, which appears in Lemma \ref{decomposition}. 

As a consequence, if the  incidence matrix of a hypergraph $\D$ is totally unimodular, then $A(\D,w)$ is standard graded for all weight functions $w$. This proves the sufficient part of Theorem \ref{main}. 
The necessary part of Theorem \ref{main} does not follow from Theorem \ref{BT} because 
of the condition $\w \in \ZZ^m$.

By Lemma \ref{decomposition},  Theorem \ref{main} follows from the following modification of Theorem \ref{BT} for  integral matrices with non-negative entries and for positive integral vectors $\w   \in \ZZ_+^m$.

\begin{Theorem} \label{modification}
 Let $M$ be an $m \times n$ integral matrix with non-negative entries. 
Then the polyhedron $\{\c   \in \RR^n|\ M\cdot \c   \ge \w \}$ has the decomposition property for all integral vector $\w \in \ZZ_+^m$ if and only if $M$  is totally unimodular.
\end{Theorem}

The proof of this theorem is given at the end of the next section. Following the proof of \cite[Theorem 19.4]{Sc}, we need to prepare some results on the relationship between the polyhedron $\{\c   \in \RR^n|\ M\cdot \c   \ge \w \}$ and the total unimodularity of $M$.

\section{Unimodularity versus integrality}

Recall that a matrix is called {\it unimodular} if each maximal minor equals 0, $\pm 1$
and that a rational polyhedron is called {\it integral} if all of its vertices are integral.

By a result of Veinott and Dantzig (see e.g. \cite[Theorem 19.2]{Sc}) unimodular integral matrices can be characterized by the integrality of associated polyhedra. This result can be modified for matrices with non-negative entries as follows.

\begin{Theorem} \label{VD}  
Let $M$ be an  integral $m \times n$ matrix with non-negative entries of full row rank. Then $M$ is unimodular if and only if the polyhedron $\{\c \in \NN^n|\  M\cdot\c  = \w  \}$ is integral for all integral vectors $\w   \in \ZZ_+^m$.
\end{Theorem}

\begin{proof}
If $M$ is unimodular, then $\{\c \in \NN^n|\ M\cdot\c  = \w  \}$ is integral  by the mentioned result of Veinott and Dantzig.
For the converse assume that $\{\c \in \NN^n|\  M\cdot\c  = \w  \}$ is integral for each integral vector $\w   \in \ZZ_+^m$. Let  $E$ be  an arbitrary non-singular maximal square submatrix of $M$. We have to prove that $\det E = \pm 1$. Since $\det E $ is an integer, $\det E^{-1} = (\det E)^{-1}$ is not an integer if $\det E \neq  \pm 1$.  On the other hand, using the Laplace expansion we can find an integral vector $\a$ such that $(\det E)^{-1}$ is a component of $E^{-1}\cdot\a$. Therefore, it suffices to show that $E^{-1}\cdot\a$ is an integral vector. \par
First, we can find an integral vector $\a' \in \ZZ_+^n$ such that $E^{-1}\cdot\a + \a' \in \ZZ_+^n$. Let $\b := E^{-1}\cdot\a + \a'$ and  $\w   := E\cdot\b$. Then $\w \in \ZZ_+^m$ because  $\b \in \ZZ_+^n$, $E$ has non-negative entries and no row of $E$ is zero. Let $\b^*$ be the vector obtained from $\b$ by adding zero-components so that  $M\cdot\b^* = E\cdot\b  = \w  $. Then $\b^* \in \{\c \in \NN^n|\ M\cdot\c   = \w  \}$ 
and $\b^*$ satisfies the maximal number of linearly independent constraints. Hence $\b^*$ is a vertex of the polyhedron $ \{\c \in \NN^n|\ M\cdot\c   = \w  \}$. By the assumption, $\b^*$ is integral. Therefore, $\b$ and hence $E^{-1}\cdot\a = \b - \a'$ are  integral vectors.
\end{proof}

The following corollary is again a modification of a result of Hoffman and Kruskal for integral matrices with non-negative entries (see e.g. \cite[Corollary 19.2a]{Sc}).
 
\begin{Corollary} \label{HK}
An integral $m \times n$ matrix $M$ with non-negative entries is totally unimodular if and only if the polyhedron $\{\c  \in \NN^n|\ M\cdot\c  \ge \w  \}$ is 
integral for all integral vectors $\w \in \ZZ_+^m$.
\end{Corollary}

\begin{proof}
Let $I$ denote the unit matrix of rank $m$. It is well-known that $M$ is totally unimodular if and only if the composed matrix $(I,M)$ is unimodular. 
On the other hand,  the vertices of the polyhedron $\{\c  \in \NN^n|\ M\cdot\c  \ge \w  \}$ are integral if and only if the vertices of the polyhedron $\{\b  \in \NN^{m+n}|\ (I,M)\cdot\b  = \w  \}$ are integral. Therefore, the assertion follows from Theorem \ref{VD}. 
\end{proof}

We shall use Corollary \ref{HK} to prove Theorem \ref{modification}. As observed in the preceding section,   Theorem \ref{main} follows from Theorem \ref{modification}.
\medskip

\noindent{\it Proof of Theorem \ref{modification}.} 
The necessity already follows from Theorem \ref{BT}.
To prove the sufficiency we assume that the polyhedron $Q = \{\c   \in \RR^n|\ M\cdot \c   \ge \w \}$ has the decomposition property for all integral vectors $\w \in \ZZ_+^m$. By Corollary \ref{HK}, we only need to show that $Q$ is integral.
Let $\a$ be an arbitrary vertex of  $Q$. Suppose that $\a$ is not integral. Let $k$ be the least common multiple of the denominators occuring in $\a$. Put $\c := k\a$.
Then  $\c$ is an integral vector in $kQ$. Therefore,  there are integral vectors $\c_1,...,\c_k\in Q$ such that $\c = \c_1 + \cdots + \c_k$. From this it follows that $\a = (\c_1 + \cdots + \c_k)/k$. Since $\a$ is a vertex of $Q$, we must have $\c_1 = \cdots = \c_k$ and hence $\a = \c_1 \in \NN^n$. \qed

\section{Remarks}

For a vector $\w$ we will use the notation $\w \gg 0$ if all components of $\w$ are large enough.
In the proof of Theorem \ref{VD} we may choose $\b \gg 0$, which implies $\w \gg 0$.
Therefore, we obtain the following stronger statement for the converse of Theorem \ref{VD}.

\begin{Theorem}
Let $M$ be an  integral $m \times n$ matrix with non-negative entries of full row rank. Then $M$ is unimodular if  the polyhedron $\{\c \in \NN^n|\  M\cdot\c  = \w  \}$ is integral for all integral vectors $\w \gg 0$   of $\ZZ_+^m$.
\end{Theorem}

Similarly as above, this result yields the following improvement of  Theorem \ref{main}, where  $w \gg 0$ means $\w \gg
0$.
 
\begin{Theorem}  
 If $A(\D,w)$ is a standard graded algebra for all weight functions $w \gg 0$, then $\D$ is an unimodular hypergraph.
\end{Theorem}

The above theorem doesn't hold for a sequence of weight functions $w \gg 0$.
In fact, if $\D$ is an unimodular hypergraph then $A(\D)$ is standard graded over $R$ by Theorem \ref{main}.
But $A(\D)$ needs not to be standard graded over $R$ if $A(\D,w)$ is standard graded over $R$ for a sequence of weight functions $w \gg 0$. This follows from the following observation.

\begin{Lemma}
For every hypergraph $\D$ there exists a number $d$ such that $A(\D,w)$ is a standard graded algebra over $R$ for all weight functions $w_F = kd$, $F \in \D$, $k \ge 1$.
\end{Lemma}

\begin{proof}
It is known that $A(\D)$ is always a finitely generated algebra over $R$  \cite[Theorem 4.1]{HHT}. Therefore, there exists a number $d$ such that the Veronese subalgebra $A(\D)^{(kd)}$ is standard graded over $R$ for all $k \ge 1$ \cite[Theorem 2.1]{HHT}. But $A(\D)^{(kd)} = A(\D,w)$ where $w$ is  the weight function $w_F = kd$, $F \in \D$. 
\end{proof}

On the other hand, one may ask whether $A(\D,w)$ is a standard graded algebra over $R$ for all weight functions $w$ if $A(\D)$ is standard graded. 

This question has a positive answer if $\D$ is a graph. In this case, if $A(\D)$ is standard graded over $R$, then $\D$ is bipartite \cite[Theorem 5.1]{HHT}, whence $A(\D,w)$ is standard graded for all weight functions $w$ by \cite[Theorem 5.4]{HHT}. However, we can give an example of a hypergraph $\D$ such that $A(\D)$ is standard graded over $R$, whereas $A(\D,w)$ is not standard graded over $R$ for all weight functions $w$.  For that we shall need the following result.
 
Let $M$ be the edge-vertex incidence matrix of a hypergraph$\D$. Let $m \times n$ be the size of $M$. One calls $\D$ a {\it Mengerian hypergraph} if
$$\min\{\a \cdot \c|\  \a \in \NN^n,\  M\cdot\a \ge \1 \}
= \max\{\b \cdot \1|\ \b \in \NN^m,\ M^T\cdot\b \le \c\}$$
for all $\c \in \NN^n$, where $\1$ denotes the vector $(1,...,1) \in \NN^m$.
One calls the hypergraph of the minimal covers of $\D$ the {\it blocker} of $\D$, which we denotes by $\D^*$. It is proved in \cite[Theorem 1.4]{HHTZ} that $A(\D)$ is standard graded over $R$ if and only if $\D^*$ is Mengerian. 

\begin{Example}\label{Mengerian 2}
{\rm Let $\D$ be the hypergraph on 5 vertices which has the edges
$$\{1,2,3\},\{1,5,6\},\{2,4,6\},\{3,4,5\}.$$
Then $\D^*$ is the simplicial complex with the edges
$$\{1,4\},\{2,5\},\{3,6\},\{1,2,3\},\{1,5,6\},\{2,4,6\},\{3,4,5\}.$$
It is known that $\D^*$ is Mengerian, whereas $\D$ is not (see e.g. \cite[Example 1.8]{HHT}).
Therefore, $A(\D)$ is standard graded over $R$. On the other hand, 
since unimodular hypergraphs are Mengerian (see e.g. \cite[Corollary 1, p. 170]{Be}), $\D$ is not unimodular. Therefore, $A(\D,w)$ is not standard graded over $R$ for all weight functions $w$.}
\end{Example}

Finally, we would like to point out that one can test totally unimodular matrices (and hence  unimodular hypergraphs) in polynomial time (see e.g. \cite[Chapter 19]{Sc}).

\end{document}